\documentclass[11pt]{article}
\usepackage[centertags]{amsmath}
\usepackage{amsfonts}
\usepackage{newlfont}
\usepackage{graphicx}
\usepackage{eurosym}
\usepackage[pdfpagelabels]{hyperref}
\usepackage{epstopdf}
\usepackage{amsmath,amssymb, xcolor}
\usepackage{caption}
\usepackage{subcaption}
\usepackage{amsfonts}
\usepackage{textcomp}
\usepackage{appendix}
\usepackage{newlfont}
\usepackage{graphicx}
\usepackage{eurosym}
\usepackage{bm}
\usepackage{mathtools}

\newcommand{\E}{\mathbf{E}}

\def\cD{{\cal D}}

\def\cF{{\cal F}}
\def\cG{{\cal G}}

\def\cP{{\cal P}}

\begin{document}

\title{{Enhanced Cauchy Schwarz inequality and some of its statistical applications}}

\author{S. Scarlatti\thanks{Dept. of Economics and Finance - University of Rome Tor Vergata, \texttt{sergio.scarlatti@uniroma2.it}}}

\maketitle

\date{}
\begin{abstract}
We present  general refinements of  the Cauchy-Schwarz and Buzano inequalities over complete inner product spaces showing these improvements have interesting statistical applications.

\medskip

\end{abstract}

\maketitle
\section{Introduction}
It is widely recognized that Cauchy-Schwarz (CS) inequality  is one of the most important inequality in mathematical analysis, if not the most important one, see (\cite{S}) and (\cite{D1}). This reputation comes from its large spectrum of applications in different fields of science ranging from physics to probability and statistics, but including also signal analysis, coding theory and computer science. It is ascertained that the Russian mathematician V. Bunyakovsky has contributed to its discovery and diffusion. It is also well known that the basic CS inequality admits several generalizations to function, matrix and operator spaces, however it is not so widely known that the basic inequality itself admits improvements. The aim of the present note is twofold:  proving  enhanced versions of the Cauchy-Schwarz and Buzano  inequalities and showing by a set of examples arising from the fields of statistics and data analysis, except for the last one, that they are easy to apply and useful to consider.\\
\section{Enhanced Cauchy-Schwarz and Buzano inequalities}
Let $(H,\langle,\rangle)$ a real or complex Hilbert space  and $||x||=\sqrt{\langle x,x\rangle}$ the vector norm, the celebrated Cauchy-Schwarz inequality states that
\begin{equation}\label{CS}
|\langle x,y \rangle|\leq ||x||\cdot ||y||\;\;\;\forall x,y\in H.
\end{equation}
In what follows only real Hilbert spaces will be contemplated. In a finite dimensional setting let us take $x,y\in H=\mathbb{R}^n$, endowed with the Euclidean norm, and consider the anti-symmetric square matrix $C=C(x,y)$ having elements defined as $(c_{ij})\equiv \frac1{\sqrt{2}}(x_iy_j-x_jy_i)$,\;$i,j=1,\ldots,n$, with $||C||_2\equiv \sqrt{\sum_{i,j}c^2_{ij}}$ its $L^2$-norm. Then CS inequality stems from the evaluation of the matrix norm, that is:
\begin{equation}\label{Frechet}
||C||_2^2=||x||^2||y||^2-\langle x,y\rangle ^2,
\end{equation}
an identity which can be traced back to Lagrange, see (\cite{S}) and (\cite{MN}).  The following proposition, which  sharpens CS inequality, is a consequence of the existence of orthogonal projections on closed subspaces of a Hilbert space, a classical theorem, see (\cite{B}). We state it in the following useful form:\\

{\bf{Proposition 2.1}}: Let $(H, \langle \cdot,\cdot\rangle )$ a real Hilbert space, $V$ any closed subspace of $H$. Let $P\equiv P_V$ denote the orthogonal projection on $V$,  that is $PH=V$,  and define the following symmetric function: 
\begin{equation}\label{D-function}
D(x,y|P)\equiv ||Px||\cdot||Py||+||x-Px||\cdot||y-Py||\;\;\forall x,y\in H.
\end{equation}
Then it holds
\begin{equation}\label{proj}
|\langle x,y\rangle |\leq D(x,y|P) \leq  ||x||\cdot||y||\;\;\forall x,y\in H.
\end{equation}
Furthermore, let  $\cP$ denote the family of all orthogonal projections relative to all possible closed subspaces of $H$, then for each pair $(x,y)$ we have
\begin{equation}\label{D^+}
D^+(x,y)\equiv \sup_{P\in \cP}D(x,y|P)=||x||\cdot ||y||
\end{equation}
and
\begin{equation}\label{D^-}
D^{-}(x,y)\equiv \inf_{P\in \cP}D(x,y|P)=|\langle x,y\rangle |.
\end{equation}
The enhanced CS inequality is therefore the statement that for any orthogonal projection $P$ it holds
\begin{equation}
D^{-}(x,y)\leq D(x,y|P)\leq D^{+}(x,y),\;\;\forall x,y \in H.
\end{equation}
{\bf{proof}}: Writing $x=Px+(I-P)x$ by orthogonality we have $||x||^2=p(x)^2+q(x)^2$,
where $p(x)\equiv ||Px||$ and $q(x)\equiv \sqrt{||x||^2-p(x)^2}=||(I-P)x||$; therefore
$$
|\langle x,y\rangle |=|\langle Px,Py\rangle +\langle (I-P)x,(I-P)y\rangle |\leq |\langle Px,Py\rangle |+|\langle (I-P)x,(I-P)y\rangle |
$$
$$
\leq p(x)p(y)+q(x)q(y) \leq \sqrt{p(x)^2+q(x)^2}\sqrt{p(y)^2+q(y)^2}= ||x||\cdot||y||
$$
by applying orthogonality, triangle inequality, the Cauchy-Schwarz inequality and noticing 
\begin{equation}\label{squaring}
(p(x)p(y)+q(x)q(y))^2+(p(x)q(y)-p(y)q(x))^2=(p(x)^2+q(x)^2)(p(y)^2+q(y)^2).
\end{equation}
The second part of the statement follows by choosing $V=\{0\}$ or $V=H$, that is $PH=\{0\}$ or $PH=H$ in order to get $D^+(x,y)$.
On the other hand the choice of $V=\{\lambda x\}$, that is of $P=P_x$, with $P_xz\equiv \langle {\frac{x}{||x||}},z\rangle \frac{x}{||x||}$ leads to $D^-(x,y)$; indeed
$$
D(x,y|P_x)= ||P_xx||\cdot||P_xy||+||x-P_xx||\cdot||y-P_xy||
$$
$$
= ||x||\cdot||P_xy||=||x||\cdot |\langle {\frac{x}{||x||}},y\rangle |=|\langle x,y\rangle |.\;\;\;\;\;\square
$$
It is worth noticing the following basic properties of the decoupling symmetric function $D(\cdot,\cdot|P)$ which has been previously introduced (we set $P^{\perp}\equiv I-P$ to shorten the notation):\\
\\ \ 
(i)$D(x,x|P)=||x||^2$ \text{and} $D(x,y|P^{\perp})=D(x,y|P)\geq 0$,\\
(ii)$D(\lambda x,\mu y|P)=|\lambda\mu| D(x,y|P)$\; $\forall \lambda , \mu \in \mathbb{R}$ \;\text{and} $D(x+x',y|P)\leq D(x,y|P)+ D(x',y|P)$,\\
(iii)$D(P^\#x,y|P)=||P^\#x||\;||P^\#y||$ for $P^\#\in\{P,P^{\perp}\}$,\\
(iv)$D((P-P^{\perp})x,y|P)=D(x,y|P)$,\\
(v)$D(x,y|P)= 0 \;\text{for all} \;(x,y)\in V\times V^{\perp}$, with $V=PH$.
\\ \ \\
{{\bf{Corollary 2.1}}: Let $(X, \langle \cdot,\cdot\rangle )$ a real Hilbert space, then for any orthogonal projection $P$ on a closed subspace it holds
\begin{equation}\label{genBuz}
2|\langle Px,y\rangle |\leq D(x,y|P)+|\langle x,y\rangle |\;\;\forall x,y\in X, 
\end{equation}
or equivalently
\begin{equation}\label{genBuzbis}
|\langle Px,y\rangle |\leq D(P^{\perp}x,y|P)+|\langle x,y\rangle |\;\;\forall x,y\in X.
\end{equation}
{{\bf{Remark 2.1}}: The less stringent inequality $2|\langle Px,y\rangle |\leq ||x||\;||y||+|\langle x,y\rangle|$ is proven in (\cite{D2}).\\
{{\bf{proof}}:} Similarly to (\cite{FK}) the basic observation is the following one:
$$
2|\langle Px,y\rangle |-|\langle x,y\rangle |\leq |\langle (P-P^{\perp})x,y\rangle |\leq D((P-P^{\perp})x,y|P)=D(x,y|P) \;\; \forall x,y\in X,
$$
the last equality coming from property (iv) of the $D$ function. Then inequality (\ref{genBuzbis}) follows from property (iii) and the very definition of the $D$ function.\;\;$\square$\\
\\ \ 
Recall that Buzano  inequality is the statement: 
\begin{equation}\label{buz}
2 |\langle x,z\rangle \langle z,y\rangle |\leq ||z||^2(||x|| ||y||+|\langle x,y\rangle|), \forall x,y,z\in X,
\end{equation}
see e.g. \cite{S}. Notice that for $x=y$ it reduces to the CS inequality. However we have:
\\ \ \\
{{\bf{Corollary 2.3}}:  Let $(X, \langle \cdot,\cdot\rangle )$ a real Hilbert space, then for all $x,y, z$ in $X$ it holds:
\begin{equation}\label{Buz}
|\langle x,z\rangle \langle z,y\rangle |\leq \sqrt{||z||^2||x||^2-\langle x,z\rangle^2}\;\sqrt{||z||^2||y||^2-\langle y,z\rangle^2}+||z||^2|\langle x,y\rangle| \; \;.
\end{equation}
{{\bf{proof}}:}W.l.o.g. we may take $z\neq 0$ and set $v\equiv \frac{z}{||z||}$. By choosing $P=P_v$ in (\ref{genBuzbis}) we obtain
\begin{equation}\label{Bu-impr}
|\langle x,v\rangle \langle v,y\rangle |\leq ||P^{\perp}_vx||\;||P^{\perp}_vy||+|\langle x,y\rangle|. 
\end{equation}
and setting $v=\frac{z}{||z||}$ in this last formula we obtain (\ref{Buz}).$\square$\\
Since in (\ref{Bu-impr}) the quantity $D(x,y|P_v)$ is not greater than the quantity $||x|| ||y||$ we have obtained an improvement w.r.t. (\ref{buz}); we refer to (\cite{LS}) for refinements of different type.\\ To grasp the content of the preceding results let us consider the following simple but instructive example:\\
{\bf{Example2.1}} Let $H=\mathbb{R}^n$, endowed with the Euclidean distance $d_n$; for $k\geq 1$ consider the subspace 
$V_{1:k}\equiv \{z|z=(z_1,z_2,\ldots,z_k,0,\ldots,0)\}$ and its orthogonal complement $V^{\perp}_{1:k}=V_{k+1:n}\equiv \{z|z=(0,\ldots,0,z_{k+1},z_{k+2},\ldots,z_n)\}$ . Define the orthogonal projection $P_kH\equiv V_{1:k}$, then we have
\begin{equation}\label{project}
D(x,y|P_k)=\sqrt{x_1^2+\ldots+x^2_k}\sqrt{y_1^2+\ldots+y_k^2}+\sqrt{x_{k+1}^2+\ldots+x^2_n}\sqrt{y_{k+1}^2+\ldots+y_n^2}
\end{equation}
and by Proposition 2.1 the following inequalities
\begin{equation}\label{nD}
|\sum_{i=1}^nx_iy_i|\leq D(x,y|P_k)\leq \sqrt{x_1^2+\ldots + x_n^2}\cdot \sqrt{y_1^2+\ldots y^2_n}\;\;\; \forall x, y \in \mathbb{R}^n
\end{equation}
hold, with matrix form of $P_k$ w.r.t the canonical basis of $\mathbb{R}^n$ given by
 $diag(I_{k\times k},{\boldsymbol{0}}_{(n-k)\times (n-k)})$. \\
The interpretation of (\ref{nD}) is straightforward: pick two arbitrary points A and B in the $n$-dimensional euclidean space and orthogonally project them on $V_{1:k}$ getting points $\tilde A$ and $\tilde B$ and on $V^{\perp}_{1:k}$ getting points $\hat A$ and $\hat B$. Then the number $d_n(0,\tilde A)d_n(0,\tilde B)+d_n(0,\hat A)d_n(0,\hat B)$ is not greater than the number $d_n(0,A)d_n(0,B)$ and not smaller than the value of the scalar product of the two vectors pointing from the origin to $A$ and $B$. It is easy to see that there are cases for which the improvement can be relevant.

\section{Applications of the enhanced Cauchy-Schwarz and Buzano inequalities}
In this section, by means of a series of different  examples, we show that sharpened CS and Buzano inequalities may find applications in applied statistics, data anlysis and numerical functional analysis. \\ \ \\
{\bf{Example3.1 (Sample covariance)}} Let $H=\mathbb{R}^n$, $e=(\frac1{\sqrt{n}},\ldots,\frac1{\sqrt{n}})$ and $u=\sqrt{n}e$.  Define $Px\equiv \langle x,e\rangle e$ for all $x\in H$, by Proposition 2.1 we have
$$
|\sum_{i=1}^nx_iy_i|\leq |\langle x,e\rangle |\cdot |\langle y,e\rangle |+\sqrt{||x||^2-|\langle x,e\rangle |^2}\cdot\sqrt{||y||^2-|\langle y,e\rangle |^2}
$$
$$
=\frac{|\langle x,u\rangle |}{\sqrt{n}}\cdot \frac{|\langle y,u\rangle |}{\sqrt{n}}+\sqrt{||x||^2-\biggl(\frac{|\langle x,u\rangle |}{\sqrt{n}}\biggr)^2}\cdot\sqrt{||y||^2-\biggl(\frac{|\langle y,u\rangle |}{\sqrt{n}}\biggr)^2}.
$$
therefore by setting $\bar x=\frac{x_1+\ldots+x_n}{n}$ we obtain the estimate:
\begin{equation}\label{product}
|\sum_{i=1}^nx_iy_i|\leq n |\bar x \bar y|+ \sqrt{||x||^2-n |\bar x|^2}\cdot\sqrt{||y||^2-n |\bar y|^2}.
\end{equation}
Notice that  by applying (\ref{product}) to $x'= x-\bar x u$ and $y'= y- \bar y u$ we get the  sample covariance estimate
\begin{equation}\label{conf}
\big|\langle x,y\rangle  -n\;\bar x \bar y\big|\leq ||x'||\;||y'||=n\cdot {\overline{s}(x)}{\overline{s}(y)},
\end{equation}
where $\overline{s}(x)^2\equiv \frac{1}{n}||x||^2-\overline x^2$. 
Furthermore, squaring (\ref{product}) gives 
\begin{equation}\label{Ovar}
(\sum_{i=1}^nx_iy_i)^2 \leq \biggl(\; n |\bar x \bar y|+ \sqrt{||x||^2-n |\bar x|^2}\cdot\sqrt{||y||^2-n |\bar y|^2}\;\biggr)^2=
\end{equation}
$$
=||x||^2||y||^2-n\biggr(|\bar x|\sqrt{||y||^2-n |\bar y|^2}-|\bar y|\sqrt{||x||^2-n |\bar x|^2}\biggl)^2.
$$
Finally we notice that, in the same vein, one may choose $H=L^2(S,\cF,Q)$, $\cF$ a sigma-algebra, $Q$ a probability measure, and $\langle X,Y\rangle \equiv \E (XY)$. For any $X\in H$ define the one-dimensional projection $PX\equiv \langle X,1\rangle 1=(\E X)1$. In this case $D(X,Y|P)=|\E X||\E Y|+\sigma_X\sigma_Y$ and (\ref{squaring}) implies
\begin{equation}\label{seeH}
(\E(XY))^2\leq ||X||^2||Y||^2-\biggr(|\E X|\sigma_Y-|\E Y|\sigma_X\biggl)^2\leq ||X||^2||Y||^2 \;\;\forall X,Y \in H,
\end{equation}
which is the main result presented  in (\cite{W}), jointly with a pair of nice applications. It must be remarked that in case both $X$ and $Y$ have zero mean the estimate (\ref{seeH})  offers no improvement over CS,   and the same it holds for centered data in the previously displayed sample estimates. This is due to the fact that the chosen $P$ is projecting onto a subspace orthogonal to 
$W=\{X\in L^2:\E(X)=0\}$.\\ \ \\
{\bf{Example3.2(Sample cross-covariance)}} Let us consider the sample cross-covariance function $\overline R_{(x,y)}(h)$ among two temporal sequences of observations $(x_1,x_2,\ldots,x_n)$ and $(y_1,y_2,\ldots,y_n)$:
$$
\overline R_{(x,y)}(h)\equiv \frac{1}{n}\sum_{t=1}^{n-h}(x_{t}-\bar x)(y_{t+h}-\bar y),
$$
with $h=1,\ldots,n-1$. In the following we shall suppose the observed data being generated by two mean-zero ergodic processes so, by assuming $n$ sufficiently large, we may take $\bar x =\bar y= 0$.  For a fixed $h$ we take $H=\mathbb{R}^{n-h}$, and choose the projection $P_k$, $k\leq n-h$, as in Example 2.1; by using the first inequality in (\ref{nD}), we have
\begin{equation}\label{cross1}
|\overline R_{(x,y)}(h)|=\frac{1}{n}|x_1\underbrace{y_{1+h}}_{z_1}+x_2\underbrace{y_{2+h}}_{z_2}+\ldots+x_{n-h}\underbrace{y_{n}}_{z_{n-h}}|\leq \frac1{n}D(x,z|P_k),
\end{equation}
with
$$
D(x,z|P_k)=\sqrt{x_1^2+\ldots+x^2_k}\sqrt{y_{1+h}^2+\ldots+y_{k+h}^2}+\sqrt{x_{k+1}^2+\ldots+x^2_{n-h}}\sqrt{y_{k+h+1}^2+\ldots+y_{n}^2}.
$$
For instance, setting $z_{r:s}\equiv (z_r,\ldots,z_s)\in \mathbb{R}^{s-r+1}$, $1\leq r\leq s\leq n$ and
 choosing in (\ref{cross1}) $k=h$\;(for $h=1,\ldots,[n/2]$),the dimension of the subspace equalizing the number of the forward time-shift steps, leads to
$$
|\overline R_{(x,y)}(h)|\leq (\frac{h}{n})\sqrt{\overline{x^2}_{1:h}}\sqrt{\overline{y^2}_{h+1:2h}}+(1-\frac{2h}{n})\sqrt{\overline{x^2}_{h+1:n-h}}\sqrt{\overline{y^2}_{2h+1:n}}\;,
$$
which is a linear combination of products of sample standard deviations on different time intervals improving over the classical CS inequality estimate.\\ \ \\
{\bf{Example3.3 (Correlations)}} We take $H=L^2(S,\cF,Q)$ as in Example 3.1. By setting $X'\equiv X-\E X$ and $Y'\equiv Y-\E Y$, CS inequality implies 
 \begin{equation}\label{cov1}
 |cov(X,Y)|=|\E(X'Y')|\leq ||X'||\;||Y'||=\sigma_X\sigma_Y,
 \end{equation}
 from which the classical bound $|\rho_{X,Y}|\leq 1$.
Consider now any sub-sigma algebra $\cG\subseteq \cF$ and the closed subspace of $H$ given by $V_{\cG}\equiv L^2(\cG,Q)$. Let $P_{\cG}H\equiv V_{\cG}$ the associated orthogonal projection, therefore  we have $P_{\cG}X=\E(X|\cG)$ for all $X\in H$. Hence $P_{\cG}X'=\E(X|\cG)-\E X$ and
$$
D(X',Y'|P_{\cG})=||\E(X'|\cG)||\cdot ||\E(Y'|\cG)||+||X'-\E(X'|\cG)||\cdot ||Y'-\E(Y'|\cG)||
$$
$$
=||\E(X|\cG)-\E X||\cdot ||\E(Y|\cG)-\E Y||+||X-\E(X|\cG)||\cdot ||Y-\E(Y|\cG)||
$$
$$
=\sigma_{\E(X|\cG)}\;\sigma_{\E(Y|\cG)}+\sqrt{\sigma^2_X-\sigma^2_{\E(X|\cG)}}\;\sqrt{\sigma^2_Y-\sigma^2_{\E(Y|\cG)}}.
$$
By Proposition 2.1 we have $|\E(X'Y')|\leq D(X',Y'|P_{\cG})\leq ||X'||\cdot||Y'||$ that is
\begin{equation}\label{cov}
|cov(X,Y)|\leq \sigma_{\E(X|\cG)}\;\sigma_{\E(Y|\cG)}+\sqrt{\sigma^2_X-\sigma^2_{\E(X|\cG)}}\;\sqrt{\sigma^2_Y-\sigma^2_{\E(Y|\cG)}} \leq \sigma_X\sigma_Y,
\end{equation}
For instance, in case $X \in \cG$ and $Y \in \cG$ then (\ref{cov}) reduces to (\ref{cov1}), while in case $X \notin \cG$ and $Y \in \cG$ then (\ref{cov}) takes the form:
\begin{equation}\label{cov3}
 |cov(X,Y)|\leq \sigma_{\E(X|\cG)}\cdot \sigma_Y\leq \sigma_X\sigma_Y.
\end{equation}
Notice that when $\cG\equiv \sigma(Y)$  and the variables $(X,Y)$  are jointly normal we have $\E(X|\cG)=\mu_X+\rho_{X,Y}\frac{\sigma_X}{\sigma_Y}(Y-\mu_Y)$, so the lower inequality in (\ref{cov3}) becomes an equality.\\
Finally,  the above estimates suggest the consideration of correlation coefficients defined as 
\begin{equation}\label{corr_G}
\rho_{X,Y}^{\cG}\equiv \frac{E(X'Y')}{D(X',Y'|P_{\cG})},
\end{equation}
(in case $\rho_{X,Y}^{\cG}=\frac{0}{0}$ its value is set to zero) which values depend on the sigma-algebra $\cG$ inherent to the problem under investigation. Indeed it holds $|\rho_{X,Y}|\leq |\rho_{X,Y}^{\cG}|\leq 1$, therefore encoding information might significantly enlarge the range of the estimated correlation. When $\cG$ is the trivial sigma-algebra the two numbers $\rho_{X,Y}$ and $\rho_{X,Y}^{\cG}$ coincide, and the random variable $\rho_{X,Y|\cG}$, defined in the standard way, also equalizes this value.\\
{\bf{Remark3.1}}:Let $X,Y \in L^2(S,\cF,Q)$, generalizing the idea of the previous example, we can use the lower inequality in (\ref{proj}) to define a P-correlation coefficient $\rho_{X,Y}(P)$ as 
\begin{equation}\label{proj_corr}
\rho_{X,Y}(P)\equiv \frac{cov(X,Y)}{D(X,Y|P)}\in[-1,1],
\end{equation}
Notice that $\rho_{X,Y}(I)=\rho_{X,Y}$ but in general $\rho_{X,Y}(P)\neq \rho_{PX,PY}$; obviously $\rho_{X,Y}(P_{\cG})=\rho_{X,Y}^{\cG}$.\\

{\bf{Example3.4 (Density divergence estimation)}}  Consider the Hilbert space  $L^2\equiv L^2(\mathbb{R}^d,dx)$ and $\text{Dens}^2\equiv L^2\cap \text{Dens}(\mathbb{R}^d)$, where $\text{Dens}(\mathbb{R}^d)$ is set of all probability density functions over $\mathbb{R}^d$. It is then possible to introduce the functional $\cD iv(\cdot,\cdot|P):\text{Dens}^2\times \text{Dens}^2\to [0,\infty]$ defined as follows
\begin{equation}\label{div}
\cD iv(f,g|P)\equiv -log\bigr(\frac{\int_{\mathbb{R}^d} f(x)g(x)dx}{D(f,g|P)}\bigl).
\end{equation}
Recalling that $D(f,f|P)=||f||^2$ we have $\cD iv(f,g|P)=0$ for $f=g$.
For $P=I$ the functional is known as the Cauchy-Schwarz divergence and has found application in density based clustering and machine learning, (\cite{JPEE}), (\cite{KLBLS}). CS divergence does not verify the triangle inequality and therefore does not define a true distance but only a pseudo-distance, however it is symmetric, positive and null on the diagonal. The CS P-divergence  introduced by (\ref{div}) shares exactly the same properties, moreover the following holds
$$
\cD iv(f,g|I)-\cD iv(f,g|P)=-log\biggl(\frac{D(f,g|P)}{||f||\;||g||} \biggr)\geq 0.
$$
Therefore, using (\ref{proj}), we have $\cD iv(f,g|I)\geq \cD iv(f,g|P)\geq 0$. Furthermore, if $(P_N)_{N\geq 1}$ is such that $||P_N u||\to ||u||$ $\forall u\in L^2$, then $D(f,g|P_N)\to ||f||\cdot ||g||$ and the CS $P_N$-divergences converge to the CS divergence. This last property allows for a  nonparametric estimation of the CS divergence among two densities based on projection estimators, see (\cite{T}), which we are going to outline fixing $d=1$. Let $(e_k(\cdot))_{k\geq 1}$ be an orthonormal basis, for any $u\in L^2$ set $u_k\equiv \langle u,e_k\rangle $, so that $||u||^2=\sum_{k=1}^{\infty}u_k^2$. For $N\geq 1$ define orthogonal projections $P_Nu\equiv \sum_{k=1}^{N}u_ke_k$, clearly $||P_N u||\to ||u||$ $\forall u\in L^2$. Suppose $(X_1,\ldots,X_n)$ be i.i.d. random variables having density $f$ and  $(Y_1,\ldots,Y_n)$ be i.i.d. random variables having density $g$ and define unbiased estimators for $f_k$ and $g_k$ respectively as $\hat f_{k,n}\equiv n^{-1}\sum_{i=1}^ne_k(X_i)$ and $\hat g_{k,n}\equiv n^{-1}\sum_{j=1}^ne_k(Y_j)$. It follows that $t_N(f)\equiv ||P_Nf||^2$ can be estimated by $\hat t_{N,n}(f)\equiv \sum_{k=1}^N \hat f_{k,n}^2$ and similarly $t_N(g)$ can be estimated by $\hat t_{N,n}(g)$. Henceforth, reconsideration of the basic Example 2.1 suggests that the CS $P_N$-divergence $T_N\equiv Div(f,g|P_N)$ may be estimated by using the statistics 
$$
\hat T_{N,n}({\bold{X}},\!{\bold{Y}})\equiv log\biggl(\frac{\sqrt{\hat t_{N,n}(f)}\sqrt{\hat t_{N,n}(g)}+\sqrt{\hat r_{N,n}(f)}\sqrt{\hat r_{N,n}(g)}}{\sum_{k=1}^{2N}\hat f_{k,n}\hat g_{k,n}}\biggr),
$$
where $\hat r_{N,n}(f)\equiv  \sum_{k=N+1}^{2N} \hat f_{k,n}^2$, so that for $n$ and $N$ sufficiently large this leads to an estimation of the CS divergence of the two densities.\\ \ \\
{\bf{Example3.5 (Numerical radius and RKHS)}} Let $X$ be a Hilbert space and $T$ a continuous linear map from $X$ into itself such that $D(T)=X$. The spectral radius and the numerical radius of $T$ are respectively  the numbers $\rho(T)\equiv \{\sup |\lambda|: \lambda\in \sigma(T)\}$ and  $w(T)\equiv\{\sup|\langle Tx,x\rangle |:x\in X, ||x||=1\}$ and it holds $\rho(T)\leq w(T)$. As the spectral radius also the numerical radius turns out to be relevant in applications, for instance in the stability analysis of numerical schemes, see  (\cite{G}), (\cite{CW}), and  in the diagnostic of MCMC algorithms, see (\cite{RR}).The numerical radius defines a norm equivalent to the operator norm $||T||\equiv \{\sup||Tx||: x\in X, ||x||=1\} $, by means of the estimates $2^{-1}||T||\leq w(T)\leq ||T||$. An upper bound for $w(T)$, better of than the previous one, is due to (\cite{K}):
\begin{equation}
w(T)\leq \frac{1}{2}(||T||+||T^2||^{1/2}),
\end{equation}
moreover it also well known that $w(T^n)\leq w(T)^n$ for all $n\geq 1$, we refer the reader to  (\cite{BDMP}) for more results on the topic. 
Let us also introduce the number $c(T)\equiv \{\inf |\langle Tx,x\rangle|: x\in X, ||x||=1\}$, known as Crawford number, hence $0\leq c(T)\leq w(T)$. We claim that, as a consequence of Corollary 2.2., the following inequality holds:
\begin{equation}\label{gap}
w(T)^2 -w(T^2)\leq ||T||^2-c(T)^2.
\end{equation}
Indeed by setting in (\ref{Buz}}) $x=Tz$ , $y=T^*z$ and taking the supremum over $S=\{z\in X: ||z||=1\}$, we obtain
$$
w(T)^2\leq \sqrt{\sup_{S}||Tz||^2-(\inf_{S} |\langle Tz,z\rangle|)^2}\;\sqrt{\sup_{S}||T^*z||^2-(\inf_{S} |\langle T^*z,z\rangle|)^2}\;+ w(T^2)
$$
$$
\leq ||T||^2-c(T)^2+w(T^2).
$$
where we have used $||T||=||T^*||$ and $c(T)=c(T^*)$. Henceforth  the estimate (\ref{gap}) bounds  the gap $\Delta^{(2)}\equiv w(T)^2-w(T^2)$ in terms of the difference between the operator norm  and the Crawford number. We conjecture that (\ref{gap}) might hold for any $n\geq 1$. The previous analysis can be extended to the context of reproducing kernel Hilbert spaces (RKHS), a central topic in Statistical Learning. To this aim let $X= H(\Omega)$, a  Hilbert space of functions defined on a non empty set $\Omega$ with a reproducing kernel $k_{\lambda}\equiv k(\cdot,\lambda)\in X$ for any $\lambda \in \Omega$. Once again consider a bounded linear operator $T:X \to X$ and unit vectors in $X$ given by $\hat k_{\lambda}=\frac{k_{\lambda}}{||k_{\lambda}||}$. The operator $T$ induces a function $\widetilde T$ on $\Omega$ defined as $\widetilde T(\lambda)\equiv  \langle T\hat k_{\lambda},\hat k_{\lambda}\rangle$, known as the Berezin transform of $T$, see (\cite{Z}) for more on the subject. The Berezin number of $T$ is then the value ber$(T)\equiv\sup\{|\widetilde T(\lambda)| : \lambda\in\Omega\}$, clearly it holds $0\leq ber(T)\leq w(T)$. Squaring ber$(T)$ not always produces a value larger than ber$(T^2)$, see e.g. (\cite{MMM}), however we obtain
\begin{equation}\label{ber}
|ber(T)^2 -ber(T^2)|\leq ||T||^2-c(T)^2.
\end{equation}
Indeed a way of reasoning similar to the one leading to (\ref{gap}) proves the estimate.\\
{\textbf{Declaration:}}This research did not receive any specific grant from funding agencies in the public, commercial, or not-for-profit sectors.

\end{document}